%% file: paper.tex
\documentclass{IEEEtran}

\usepackage{cite}

\usepackage{amsthm,amsmath,amssymb,amsfonts}
\usepackage{bm}
\usepackage{algorithmic}
\usepackage[ruled, noline]{algorithm2e}
\usepackage{graphicx}

\usepackage{array}
\usepackage{multirow}
\usepackage{diagbox}

\usepackage{subcaption}
\usepackage{xcolor}

\usepackage{textcomp}
\def\BibTeX{{\rm B\kern-.05em{\sc i\kern-.025em b}\kern-.08em
T\kern-.1667em\lower.7ex\hbox{E}\kern-.125emX}}

\usepackage{lipsum}

\theoremstyle{plain}
\newtheorem{theorem}{Theorem}

\newtheorem{corollary}{Corollary}

\newtheorem{remark}{Remark}

\theoremstyle{definition}

\def\({\left(}
\def\){\right)}
\def\[{\left[}
\def\]{\right]}

\def\dbf{{\bf d}}

\def\Kbf{{\bf K}}  
\def\Lbf{{\bf L}}  
\def\Qbf{{\bf Q}}  

\def\ubf{{\bf u}}
\def\xbf{{\bf x}}
\def\ybf{{\bf y}}

\def\betabf{{\bm{\beta}}}
\def\Phibf{\bm{\Phi}}
\def\Rmbb{\mathbb{R}}  
\def\Hcal{\mathcal{H}}  
\def\Rcal{\mathcal{R}}  
\def\Scal{\mathcal{S}}

\newcommand{\SFpair}{\{\Phibf_\xbf, \Phibf_\ubf\}}
\newcommand{\OFqple}{\{\Phibf_{\xbf\xbf}, \Phibf_{\ubf\xbf}, \Phibf_{\xbf\ybf}, \Phibf_{\ubf\ybf}\}}

\def\mat#1{\begin{bmatrix}#1\end{bmatrix}}
\def\t{[t]}
\def\tn{[t+1]}
\def\tm{[t-1]}
\def\alg#1{Algorithm~\ref{alg:#1}}
\def\fig#1{Fig.~\ref{fig:#1}}
\def\subfig#1#2{Fig.~\ref{fig:#1}(\subref{subfig:#1-#2})}
\def\sec#1{Section~\ref{sec:#1}}
\def\apx#1{Appendix~\ref{apx:#1}}
\def\tab#1{Table~\ref{tab:#1}}
\def\thm#1{Theorem~\ref{thm:#1}}
\def\eqn#1{\eqref{eqn:#1}}

\def\st{{\rm s.t.}}
\def\OptConsSep{&&\quad}

\newcommand{\OptMin}[2]{
\begin{alignat}{2}
\text{minimize}\ &\ #1 \nonumber \\
\st\ #2
\end{alignat}
}

\newcommand\OptCons[3]{
&\ #1
\ifx\\#2\\ \else \OptConsSep #2 \fi%
\ifx\\#3\\ \nonumber \else \label{eqn:#3} \fi%
}

\def\algnoindent{\item[$\bullet$] }

\definecolor{mblue}{rgb}{0.075,0.541,0.855}

\renewcommand{\paragraph}[1]{

\vspace*{0.25\baselineskip}\noindent{\bf #1}}

\title{Synthesis to Deployment:\\ Cyber-Physical Control Architectures}

\author{Shih-Hao~Tseng
and~James~Anderson
\thanks{
This paper was presented in part at IEEE American Control Conference, July 1--3, 2020.
}
\thanks{Shih-Hao Tseng is with Division of Engineering and Applied Science, California Institute of Technology, Pasadena, CA 91125 USA (e-mail: shtseng@caltech.edu).}
\thanks{James Anderson is with Department of Electrical Engineering, Columbia University, New York, NY 10027 USA (e-mail: james.anderson@columbia.edu).}
}

\begin{document}

\maketitle

\bstctlcite{IEEE_BSTcontrol}

\setlength{\textfloatsep}{5pt plus 2pt minus 4pt}

\input{abstract}

\input{introduction}

\input{notation}
\input{background}

\input{implementation}

\input{architecture}

\input{original}

\input{centralized}

\input{global-state}

\input{distributed}
\input{comparison}
\input{future-directions}

\appendices
\input{appendices}

\bibliographystyle{IEEEtran}
\bibliography{Test}

\end{document}

%% file: abstract.tex
\begin{abstract}
We consider the problem of how to deploy a controller to a (networked) cyber-physical system (CPS). Controlling a CPS is an involved task, and synthesizing a controller to respect sensing, actuation, and communication constraints is only part of the challenge. In addition to controller synthesis, one should also consider how the controller will be incorporated  within the CPS. Put another way, the cyber layer and its interaction with the physical layer need to be taken into account.

In this work, we aim to bridge the gap between theoretical controller synthesis and practical CPS deployment. We adopt the system level synthesis (SLS) framework to synthesize a controller, either state-feedback or output-feedback, and provide deployment architectures for the standard SLS controllers. Furthermore, we derive new controller realizations for open-loop stable systems and introduce  different state-feedback and output-feedback architectures for deployment, ranging from fully centralized to fully distributed. Finally, we compare the trade-offs among them in terms of robustness, memory, computation, and communication overhead.
\end{abstract}

%% file: introduction.tex
\section{Introduction}\label{sec:introduction}

\begin{figure*}
\centering
\includegraphics{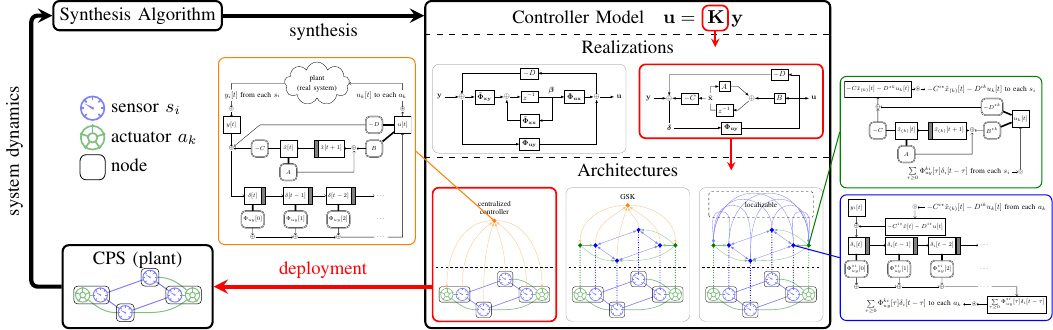}
\caption{A model-based system control scheme consists of two phases--synthesis and deployment.  Here we employ SLS in the synthesis phase to obtain a controller. The focus of this work is how to deploy $\Kbf$ to the underlying CPS.
}
\label{fig:flow-chart}
\vspace*{-\baselineskip}
\end{figure*}

\IEEEPARstart{W}{e} consider a linear time-invariant (LTI) system with a set of sensors $s_i, i = 1, \dots, N_y$, and a set of actuators $a_k, k = 1, \dots, N_u$ distributed across a network, that eveolves according to the dynamics
\begin{subequations}\label{eqn:system-dynamics}
\begin{align}
x\tn =&\ A x\t + B u\t + d_x\t \label{eqn:system-dynamics-x} \\
y\t =&\ C x\t + D u\t + d_y\t \label{eqn:system-dynamics-y}
\end{align}
\end{subequations}
where $x[t] \in \Rmbb^{N_x}$ is the state at time $t$, $y[t] \in \Rmbb^{N_y}$ is the measured output at time $t$, $u[t] \in \Rmbb^{N_u}$ is the control action at time $t$, and $d_x[t] \in \Rmbb^{N_x}$, $d_y[t] \in \Rmbb^{N_y}$ are the disturbances.
The goal of this paper is to address how a feedback controller (we will discuss state and output feedback) can be deployed to this networked system and what the corresponding cyber-physical structures and trade-offs are.

A model-based approach to control design involves two phases: the \emph{synthesis phase} and the \emph{deployment phase}, as illustrated in \fig{flow-chart}. The synthesis phase aims to derive the desired controller model using some synthesis algorithm (e.g. $\mu$, $\mathcal H_2$, $\mathcal H_{\infty}$, $L_1$, etc.) based on the system dynamics model~\eqn{system-dynamics}, a suitable objective function, and operating constraints on sensing, actuation, and communication capabilities.
The \emph{optimal} controller in the model-based sense is thus the controller model achieving the best objective value. In some cases the controller model is associated with a realization, for example, the block diagram representation of the ``central controller'' in $\mathcal H_2$ and $\mathcal H_{\infty}$ control \cite{zhou1996robust}.

Synthesizing an optimal controller for a cyber-physical system (CPS) is, in general, a daunting task due to the networked nature of sensors and actuators. Unlike centralized (single-entity) control, where sensors and actuators are presumably readily accessible, networked/distributed control is limited in measurement and actuation. Recently, the \emph{system level synthesis } (SLS) framework was proposed to facilitate controller synthesis for large-scale networked systems \cite{doyle2017system,wang2019system,AndDLM19}. By directly synthesizing desired closed-loop responses, SLS can easily incorporate system level constraints that reflect the system's networked nature, such as disturbance localization  \cite{WanMD18}, state, and input constraints \cite{CheA19}. More importantly, SLS synthesizes and optimizes the linear controller model while preserving the enforced closed-loop (distributed/network) structure.
The controller admits multiple (mathematically equivalent) control block diagrams and state space realizations \cite{wang2019system,anderson2017structured,ZheFPNK20}.

On the other hand, the deployment phase (often referred to as implementation) is concerned with how to map the derived controller realization to the physical system.
It is usually possible to implement one controller realization by multiple \emph{architectures}. Although all architectures lead to the ``optimal'' controller, they can differ in aspects other than the objective, for example; memory requirements, robustness to failure, scalability, and financial cost.
Compared to traditional centralized control systems, establishing the correct architecture and accounting for the suitable trade-offs is vital in the CPS setting.
Typically, control engineers pay less attention to deployment than they do synthesis. Our goal is to develop a theory for deployment, i.e., provide a systematic theory for controller implementation.

Approaches to deployment vary greatly in the literature.
In the control literature, most work gives the controller design at the realization level and implicitly relies on some interface provided by the underlying CPS for deployment \cite{fink2011robust,schwager2011eyes,AraMATJ14}. We draw a distinction between realization theory and deployment. Realization refers to  the mapping of an input-output model, i.e., a transfer matrix, to a state-space model. In the linear centralized setting realization theory is completely understood. In the distributed setting results are more sparse~\cite{VasE15,LesKR13}. While we do provide a realization result in Section~\ref{sec:realization}, our goal is to convince the control community that in order to deploy a control system we need to go beyond a mathematical realization of the controller and consider the underlying system architecture. Providing the appropriate abstraction and programming interface is itself a design challenge \cite{lee2008cyber,khaitan2014design,lee2015cyber,hu2016robust}.
Recognizing the communication limitations imposed by the cyber layer, there is also a long literature on networked control systems \cite{cloosterman2010controller, zhang2012network, zhang2015survey} and event-triggered control \cite{eqtami2010event,dimarogonas2011distributed,mazo2011decentralized,kartakis2017communication}. These studies usually take the cyber structure and properties (such as delay and packet dropout rate) as given and investigate which control policy is most effective. In this work, however, we focus on which architecture -- including both the cyber and physical structures -- should be deployed to have the best operational properties. Therefore, our work serves as a foundation for those upper level control policies to build upon.
Some papers also examine the deployment down to the circuit level \cite{shao2006new,Jer14}. For instance, a long-standing research area in control has looked at controller implementation using passive components \cite{smith2002synthesis,yuce2006ccii,casciati2009passive,chen2009missing}.
The networking/system community, on the other hand, mostly adopts a bottom-up instead of a top-down approach to system control. It usually involves some carefully designed gadgets/protocols and a coordination algorithm \cite{hill2000system,hamed2018chorus,dhekne2019trackio}. The survey paper~\cite{ChiLCD07} studies top down versus bottom up as well as horizontal decomposition architectures in networking.
The work in~\cite{YadSV10} is perhaps closes in spirit to our work, the authors observe that different architectures give rise to systems with varying performance levels.

\subsection{Contributions and Organization}
In this work we take an alternative approach to deployment, which lies between the realization and the circuit level. Rather than binding the design to some specific hardware, we specify a set of \emph{basic components}, i.e. components with well defined functionality that systems can be built from, and use these to implement the derived controller realization. As such, we can easily map our architectures to the real CPS.

The paper is organized as follows. In \sec{background}, we briefly review SLS, propose new, simpler state-feedback and output-feedback control realizations (for the open-loop unsatble setting), and show they are internally stabilizing.
We then discuss the gap between a mathematical realization and a practical implementation in \sec{implementation} by characterizing equivalent realizations, formulating the implementation/deployment requirements, and identifying a set of universal basic components for CPS deployment.
Leveraging the basic components, we propose different partitions and their corresponding deployment architectures in \sec{architecture}, namely, the centralized (\sec{centralized}), global state (\sec{global-state}), and four different distributed (\sec{distributed}) architectures.
Then, in Section \sec{comparison} we discuss the trade-offs made by the architectures on robustness, memory, computation, and communication.

%% file: notation.tex
\subsection*{Notation and Terminology}

Let $\Rcal\Hcal_{\infty}$ denote the set of stable rational proper transfer matrices, and $z^{-1}\Rcal\Hcal_{\infty} \subset \Rcal\Hcal_{\infty}$ be the subset of strictly proper, stable transfer matrices. Lower- and upper-case letters (such as $x$ and $A$) denote vectors and matrices respectively, while bold lower- and upper-case characters and symbols (such as $\ubf$ and ${\Phibf_\ubf}$) are reserved for signals and transfer matrices. Let $A^{ij}$ be the entry of $A$ at the $i^{\rm th}$ row and $j^{\rm th}$ column. We define $A^{i\star}$ as the $i^{\rm th}$ row and $A^{\star j}$ the $j^{\rm th}$ column of $A$. We use ${\Phi_u}[\tau]$ to denote the $\tau^{\rm th}$ spectral element of a transfer function ${\Phibf_\ubf}$, i.e., ${\Phibf_\ubf} = \sum\limits_{\tau=0}^{\infty} z^{-\tau} {\Phi_u}[\tau]$. The set $\left\lbrace x_i[t] \right\rbrace_i$ represents the elements $x_i[t]$ indexed by $i$.

We briefly summarize below the major terminology:
\begin{itemize}
\item {\bf Controller model:} a frequency-domain (linear) mapping from the output signal (provided by a sensor) to the control signal (destined for an actuator).
\item {\bf Realization:} a state space realization of the controller model together with a block diagram.
\item {\bf Architecture:} a cyber-physical structure of the controller built from a set of basic components (which we define in \sec{implementation-components}).
\item {\bf Node:} a physical sub-system of a CPS that can host and execute basic components.
\item {\bf Synthesize:} derive a controller model (and some realization).
\item {\bf Deploy/Implement:} map a controller model (through a realization) to an architecture.
\end{itemize}

%% file: background.tex
\section{Synthesis Phase}\label{sec:background}
We briefly review distributed control using system level synthesis, describing the ``standard'' SLS controller realizations and draw attention to their respective architectures.
New controller realizations for open-loop stable systems are derived, which in certain settings have favorable properties in the deployment stage.

\subsection{System Level Synthesis (SLS) -- An Overview}
To synthesize a closed-loop controller for the system described in \eqn{system-dynamics}, SLS introduces the \emph{system response} transfer matrices: $\SFpair$ for state-feedback and $\OFqple$ for output-feedback systems.
The system response is the closed-loop mapping from the disturbances to the state $\xbf$ and control action $\ubf$, under the feedback policy $\ubf = \Kbf \ybf$. Full derivations and proofs of material in this section can be  found in~\cite{wang2019system,AndDLM19}.

\subsubsection{Output-Feedback}
\begin{figure}
\centering
\includegraphics{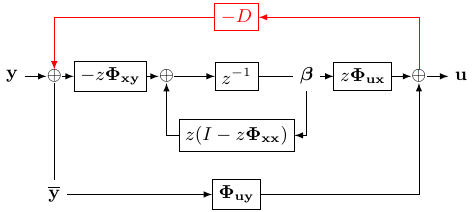}\\
\vspace*{0.5\baselineskip}
\includegraphics{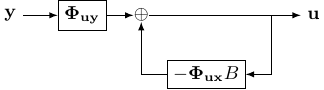}\\
\caption{Output-feedback SLS control block diagrams derived in \cite{wang2019system}. The top diagram is applicable to general systems. When  the plant is strictly proper, the path in red can be  removed. When the plant is open-loop stable \emph{and} $D = 0$, a simplified diagram (bottom) is possible. The diagrams, together with~\eqref{eqn:SLS-output-feedback-ybar}, form a realization.
}
\label{fig:SLS-output-feedback-controller-diagram}
\end{figure}

Consider the LTI system defined by~\eqn{system-dynamics}
, the system response matrices are defined by
\begin{align*}
\mat{
\xbf \\ \ubf
}
=
\mat{
\Phibf_{\xbf\xbf} & \Phibf_{\xbf\ybf} \\
\Phibf_{\ubf\xbf} & \Phibf_{\ubf\ybf}
}
\mat{
\dbf_\xbf \\ \dbf_\ybf
}
\end{align*}
where $\dbf_\xbf$ and $\dbf_\ybf$ denote the disturbance to the state and output respectively. Let $\Lbf = (\Kbf^{-1} - D)^{-1}$, the system response transfer functions take the following form:
\begin{align}
\Phibf_{\xbf\xbf} =&\ (zI-A-B \Lbf C)^{-1}, &\quad
\Phibf_{\xbf\ybf} =&\ \Phibf_{\xbf\xbf} B \Lbf  \label{eqn:OFreal} \\
\Phibf_{\ubf\xbf} =&\ \Lbf C \Phibf_{\xbf\xbf} ,&\quad
\Phibf_{\ubf\ybf} =&\ \Lbf + \Lbf C \Phibf_{\xbf\xbf} B \Lbf. \nonumber
\end{align}
An SLS problem is a convex program defined by the system response transfer functions. In it's most general setting, an SLS problem takes the form:
\OptMin{
g(\Phibf_{\xbf\xbf}, \Phibf_{\xbf\ybf},
\Phibf_{\ubf\xbf}, \Phibf_{\ubf\ybf})
}{
\OptCons{
\mat{zI-A & -B}
\mat{
\Phibf_{\xbf\xbf} & \Phibf_{\xbf\ybf} \\
\Phibf_{\ubf\xbf} & \Phibf_{\ubf\ybf}
}
=
\mat{I & 0}
}{}{of-constraint1}\\
\OptCons{
\mat{
\Phibf_{\xbf\xbf} & \Phibf_{\xbf\ybf} \\
\Phibf_{\ubf\xbf} & \Phibf_{\ubf\ybf}
}
\mat{zI-A \\ -C}
=
\mat{I \\ 0}
}{}{of-constraint2}\\
\OptCons{
\Phibf_{\xbf\xbf}, \Phibf_{\xbf\ybf}, \Phibf_{\ubf\xbf} \in \frac{1}{z}\Rcal\Hcal_{\infty},\ \
\Phibf_{\ubf\ybf} \in \Rcal\Hcal_{\infty}
}{}{}\\
\OptCons{
\mat{
\Phibf_{\xbf\xbf} & \Phibf_{\xbf\ybf} \\
\Phibf_{\ubf\xbf} & \Phibf_{\ubf\ybf}
} \in \Scal.
}{}{}
}
Observe that as as long as $g$ is convex, this is convex problem because i) constraints {\eqn{of-constraint1} and \eqn{of-constraint2}}
are affine in the decision variables, ii) the stability and strictly proper constraints are convex, and iii) $\Scal$ encodes localization constraints which are shown to be convex in~\cite{WanMD18}.\footnote{For the sake of this paper, the constraint set $\mathcal S$ does nothing other than impose sparsity constraints on the spectral components of the system response matrices. Later it will be used to impose a finite impulse response on the system response.}
The core SLS result states that the affine space parameterized by $\OFqple$ in~\eqn{of-constraint1}--\eqn{of-constraint2} parameterizes all system responses~\eqn{OFreal} achievable by an \emph{internally stabilizing} controller~\cite{wang2019system,AndDLM19}. Moreover, any feasible system response can realize an internally stabilizing controller via
\begin{align}
\Kbf =&\
\(
\( \Phibf_{\ubf\ybf} - \Phibf_{\ubf\xbf} \Phibf_{\xbf\xbf}^{-1} \Phibf_{\xbf\ybf} \)^{-1}
+ D
\)^{-1}.
\label{eqn:of-K}
\end{align}
Note that the constraint set $\Scal $ imposes sparsity constraints (amongst other things) on the system response transfer matrices. The controller defined in~\eqn{of-K} will not inherit this sparsity and will thus be unstructured -- from an implementation perspective, this is problematic. However, the realization
\begin{align}
\beta\tn =&\
- \sum_{\tau \geq 2}\Phi_{xx}[\tau]\beta[t + 2 - \tau] -
\sum_{\tau \geq 1}\Phi_{xy}[\tau]\overline{y}[t + 1 - \tau], \nonumber \\
u\t =&\
\sum_{\tau \geq 1}\Phi_{ux}[\tau]\beta[t + 1 - \tau] +
\sum_{\tau \geq 0}\Phi_{uy}[\tau]\overline{y}[t - \tau], \nonumber \\
\overline{y}\t =&\ y\t - D u\t,
\label{eqn:SLS-output-feedback-ybar}
\end{align}
corresponding to the block diagram shown in~\fig{SLS-output-feedback-controller-diagram} maintains the localized structure of the system response.

\subsubsection{State-Feedback}
\begin{figure}
\centering
\includegraphics{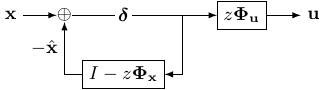}
\caption{State-feedback SLS control block diagram derived in \cite{doyle2017system} corresponding to~\eqn{SLSdefault}.  }
\label{fig:SLS-state-feedback-controller-diagram}
\end{figure}

When the state information is available to the controller, the control action is given by $\ubf = \Kbf \xbf$, and the system response is characterized by just two transfer matrices $\SFpair$, where
\begin{align*}
\Phibf_{\xbf\xbf} &\rightarrow \Phibf_\xbf := (zI-A-B \Kbf )^{-1},\\
\Phibf_{\ubf\xbf} &\rightarrow \Phibf_\ubf := \Kbf\Phibf_\xbf.
\end{align*}
The resulting SLS problem greatly simplifies; constraint~\eqn{of-constraint2} disappears, and the $2\times 2$ block transfer matrices in the remaining constraints reduce to $2\times 1$ block transfer matrices. Any feasible system response $\SFpair$ can be used to construct an internally stabilizing controller via $\Kbf = \Phibf_{\ubf}\Phibf_{\xbf}^{-1}$. As with the output-feedback case, this controller is not ideal for implementation. Instead, the controller depicted in~\fig{SLS-state-feedback-controller-diagram} inherits the structure imposed on $\SFpair$. Analogously to \eqn{SLS-output-feedback-ybar}, the controller dynamics are given by
\begin{subequations}\label{eqn:SLSdefault}
\begin{align}
\delta\t = x\t - \hat{x}\t, \quad
u\t  = \sum
\Phi_u[\tau]\delta[t + 1 - \tau],\\
\hat{x}\tn  = \sum
\Phi_x[\tau]\delta[t + 2 - \tau].
\end{align}
\end{subequations}

This paper focuses on how to actually implement and deploy such controllers into cyber-physical systems. We will characterize the necessary hardware resources, and dig further into the realizations of the  controllers and develop architectures for them.

\subsection{Simplified Controller Realizations}\label{sec:realization}

Our first result is to show that when the system to be controlled is open-loop stable, there exist alternative controller realizations that simplify the control law computation with fewer convolutions. In some circumstances, the simplified realizations may be more desirable than those described previously.

The controller realizations given by the pairs \eqn{SLS-output-feedback-ybar} and \fig{SLS-output-feedback-controller-diagram}, and \eqn{SLSdefault} and \fig{SLS-state-feedback-controller-diagram}
are preferrable to~\eqn{of-K} and $\Kbf = {\Phibf_\ubf}{\Phibf_\xbf^{-1}} $. In the output-feedback setting, four convolution operations involving  $\OFqple$ are required, and in the state-feedback setting, two convolutions involving $\SFpair$ are used. We will now show that the number of convolutions can be reduced when the plant is open-loop stable. In the sequel, we will explicitly consider the savings/trade-offs in terms of memory and computation.

\begin{figure}
\centering
\subcaptionbox{State-feedback realization corresponding to~\eqn{state-feedback-controller_time}.
\label{fig:state-feedback-controller-diagram}
}{
\includegraphics{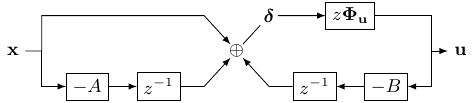}
}
\subcaptionbox{Output-feedback realization corresponding to~\eqn{output-feedback-controller_time}. The  feedback path (in red) can be removed when $D = 0$.
\label{fig:output-feedback-controller-diagram}}{
\includegraphics{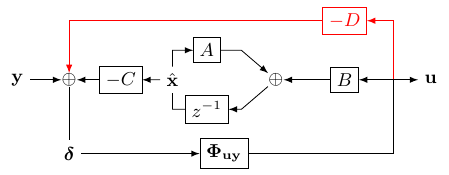}
}
\caption{SLS controller realizations corresponding to \thm{state-feedback-controller}.
}
\label{fig:simpler-realizations}
\end{figure}

\begin{theorem}\label{thm:state-feedback-controller}
Assume that  $A \in  \Rmbb^{N_x \times N_x}$ in~\eqn{system-dynamics} is Schur stable. The  following are true
\begin{itemize}
\item For the  state-feedback problem with the plant~\eqn{system-dynamics-x}, the dynamic state-feedback controller $\ubf = \Kbf \xbf$ realized via
\begin{subequations}\label{eqn:state-feedback-controller_time}
\begin{align}
\delta\t =&\ x\t - A x\tm - B u\tm,
\label{eqn:state-feedback-state-space-delta}\\
u\t =&\ \sum\limits_{\tau\geq 1} {\Phi_u}[\tau]\delta[t+1-\tau],
\label{eqn:state-feedback-state-space-u}
\end{align}
\end{subequations}
is internally stabilizing and is realized by the block diagram in \fig{state-feedback-controller-diagram}.
\item  For the  output-feedback problem with the plant~\eqn{system-dynamics}, the dynamic output-feedback controller $\ubf = \Kbf \ybf$ realized via
\begin{subequations}\label{eqn:output-feedback-controller_time}
\begin{align}
\hat{x}\tn =&\ A \hat{x}\t + B u\t,
\label{eqn:output-feedback-state-space-hatx}\\
\delta\t =&\ y\t - C \hat{x}\t - D u\t
\label{eqn:output-feedback-state-space-delta}\\
u\t =&\ \sum\limits_{\tau\geq 0} {\Phi_{uy}}[\tau]\delta[t-\tau],
\label{eqn:output-feedback-state-space-u}
\end{align}
\end{subequations}
is internally stabilizing and is described by the block diagram in \fig{output-feedback-controller-diagram}.
\end{itemize}
\end{theorem}

\begin{proof} See \apx{output-feedback-internal-stability}.
\end{proof}
In the state-feedback setting, it is possible to internally stabilize some open-loop unstable systems by the controller~\eqn{state-feedback-controller_time}. Consider a decomposition of the system matrix such that $A = A_u + A_s$ where $A_u$ is unstable and $A_s$ is Schur stable. Then using the robustness results from~\cite{MatWA17}, the controller designed for $(A_s, B)$ will stabilize $(A,B)$ if  $\|A_u \Phibf_\xbf\|<1$ for any induced norm. Of course, the ``standard'' SLS controllers do not have such an open-loop stability requirement.

\begin{remark}
The controller realizations given in Theorem~\ref{thm:state-feedback-controller}  require only one convolution operation. In contrast, the standard SLS realizations require at least two convolutions. Both realizations in \thm{state-feedback-controller} preserve the structure of the system response matrices as can be seen in ~\fig{simpler-realizations}.
\end{remark}

\input{FIR}

In the remainder of the paper we consider how a FIR controller can be implemented by a set of basic components, and compare several architectures in terms of their memory and compute requirements.

%% file: FIR.tex
\subsubsection{FIR Limitations}
In both controllers in \thm{state-feedback-controller} the summation has infinite support over the non-negative reals, i.e., they have an infinite impulse response. However, from an optimization perspective, and a performance perspective, it is beneficial to enforce a finite impulse response  (FIR) constraint which we encode in $\Scal$. The theorem above holds for both IIR and FIR controllers, however, for the remainder of the paper, we will assume all controllers are FIR with horizon length $T$. See~\cite{YuWA21} for work that removes the FIR constraint and produces a tractable convex optimization problem.  Work in~\cite{JenB20} studies continuous distributed parameter SLS controllers.

\begin{remark}\label{rem:fir}
The realization~\eqn{state-feedback-controller_time} can mix FIR and IIR blocks. When $\Phibf_\ubf$ has finite impulse response (FIR) with horizon $T$, we have
\begin{align*}
u\t =&\ \sum\limits_{\tau=1}^{T} {\Phi_u}[\tau]\delta[t+1-\tau]
\end{align*}
regardless of whether ${\Phibf_\xbf}$ is FIR or not.
\end{remark}

Although we wont's pursue this further here, we note that the horizon of the FIR filters can be optimized over. Thus,  the trade-off between hardware cost and performance cost can be seen as a pareto optimal curve. The reason this is possible is that the ceiling and floor functions
\begin{equation*}
\lceil x \rceil := \inf \{ w\in \mathbb Z~|~w\ge x \}, \quad \lfloor x \rfloor := \sup \{ w\in \mathbb Z~|~x\ge w \},
\end{equation*}
are quasilinear, i.e.,
\begin{equation*}
\lceil \theta x + (1-\theta)y \rceil = \max \{\lceil x \rceil,\lceil y \rceil \}
\end{equation*}
for all $\theta \in [0,1]$, and $x,y$ in an interval on $\Rmbb$. As such we can make  the horizon length $T$ a decision variable and include it in the synthesis objective function or include it in the constraints. Solutions to quasiconvex convex optimization problems can be obtained using a bisection algorithm to within $\epsilon$ accuracy by solving at most $\lceil \log_2(u-l)\epsilon^{-1} \rceil $ convex programs, where $p^{\star}$ (the optimal value) satisfies $l \le p^{\star} \le u$, \cite{BoyV04}. Alternatively, for problems with convex constarints and a quasiconvex objective, the subgradient method is convergent~\cite{Kiw01}.

In~\cite{YuWA21} the FIR restriction is limited in the setting of $\mathcal H_2$ optimal distributed control using SLS. While we stick with FIR filters in this paper, it should be noted that for certain choices of system norm, the tail of the impulse response can be bounded. Following~\cite{BalB92}, we consider the $\mathcal L_1$ norm, i.e., the induced $\ell_{\infty}\rightarrow \ell_{\infty}$ norm. It is enough to provide the result for multi-input, single-output systems as the MIMO result follows by taking the maximum over all output channels. Consider a minimal LTI system $G= (A,B,C,D)$ with state dimension $N_x$. Let $P$ and $Q$ be positive definite solutions to the discrete-time Lyapunov equations:
\begin{align*}
APA^T-P=-BB^T, \quad A^TQA-Q=-C^TC,
\end{align*}
corresponding to $G$. The most informative bound for the error between the IIR system $G$ and its FIR truncation of horizon length $N$ as $G_N$ is:\footnote{Note that we previously used $T$ for the FIR horizon length. However, that notation would confuse $A^T$ with the transpose of $A$.}
\begin{equation*}
\|G_N\|_{\mathcal{L}_1}+\sigma_1(Z_N)\le \|G\|_{\mathcal{L}_1} \le \|G_N\|_{\mathcal{L}_1} + 2\sum_{i=1}^{N_x}\sigma_i(Z_N),
\end{equation*}
where $\sigma_i(\cdot)$ denotes the $i^{\text{th}}$ singular value, and $Z_N = Q^{\frac{1}{2}}A^{N}P^{\frac{1}{2}}$.
Moreover, the ratio between the upper and lower bound less than or equal to $2N_x$. The upper and lower bounds can be used to determine a value for the filter horizon $N$ in order to ensure the truncation  error is less than an arbitrary constant $\epsilon$. From this we can conclude that the restriction to FIR controllers does not necessarily come at the expense of system performance.

%% file: implementation.tex
\section{Controller Implementation}\label{sec:implementation}

There is a gap between a mathematical realization and a cyber-physical implementation of the same controller. There are multiple ways to implement a device that could enforce the control law behind the mathematical realization. At the same time, each implementation requires a different amount of computation and storage resources. As a result, we face different cost trade-offs when choosing among different implementations.

When it comes to networked CPS, the situation is more puzzling. The physical structure of large networked CPSs is usually immutable in how their components interconnect due to the sunk cost or geographical restrictions. Accordingly, communication becomes an inevitable and expensive operation in time and link capacity. Therefore, in addition to computation and storage costs, we are also subject to communication constraints when we try to enforce a given mathematical realization in a networked CPS.
To highlight the unique challenges faced in CPS controller implementation, we refer to implementing a realization into a CPS as \emph{deployment}.

In the following, we characterize equivalent realizations, formulate these implementation/deployment requirements, and, accordingly, identify a set of universal basic components for CPS deployment.

\subsection{Equivalent Realizations}

Since a controller can be realized in various forms, before we dive into implementation/deployment requirements, we characterize the equivalent controller realizations in the frequency domain. A mathematical realization of the controller $\Kbf$ can be written as $\ubf = \Kbf \ybf$. In general, a controller could have some internal state $\betabf$, which is different from the state of the plant $\xbf$, and we can even express general controller realizations as
\begin{align}
\mat{\ubf\\ \betabf} =
\mat{
\Kbf_{\ubf\ybf} & 0 & \Kbf_{\ubf\betabf}\\
\Kbf_{\betabf\ybf} & \Kbf_{\betabf\ubf} & \Kbf_{\betabf\betabf}
}
\mat{\ybf\\ \ubf\\ \betabf}
\label{eqn:controller-realization-with-internal-state}
\end{align}
similar to the definitions of realizations in \cite{tseng2021realization,tsengsubgeneral}.

When $\ubf = \Kbf \ybf$ and \eqn{controller-realization-with-internal-state} are equivalent, we have the following theorem
\footnote{We can also derive \thm{K-equivalent} using the transformation technique in \cite{tseng2021realization,tsengsubgeneral}.}
:
\begin{theorem}\label{thm:K-equivalent}
Suppose the controller $\Kbf$ such that $\ubf = \Kbf \ybf$ can be implemented by the realization as in \eqn{controller-realization-with-internal-state}, then
\begin{align*}
\Kbf = (I - \Qbf\Kbf_{\betabf\ubf})^{-1}(\Kbf_{\ubf\ybf}+\Qbf\Kbf_{\betabf\ybf})
\end{align*}
where
$\Qbf = \Kbf_{\ubf\betabf}(I-\Kbf_{\betabf\betabf})^{-1}$.
\end{theorem}

\begin{proof}
From \eqn{controller-realization-with-internal-state}, we can solve for $\betabf$ by
\begin{align*}
(I-\Kbf_{\betabf\betabf})\betabf = \Kbf_{\betabf\ybf}\ybf + \Kbf_{\betabf\ubf} \ubf
\end{align*}
and hence
\begin{align*}
\betabf = (I-\Kbf_{\betabf\betabf})^{-1}(\Kbf_{\betabf\ybf}\ybf + \Kbf_{\betabf\ubf} \ubf).
\end{align*}
Substitute $\betabf$ into the equation for $\ubf$ yields
\begin{align*}
\ubf = \Kbf_{\ubf\ybf}\ybf + \Kbf_{\ubf\betabf} \betabf
= (\Kbf_{\ubf\ybf} + \Qbf \Kbf_{\betabf\ybf}) \ybf + \Qbf\Kbf_{\betabf\ubf} \ubf,
\end{align*}
which leads to
\begin{align*}
\ubf = (I - \Qbf\Kbf_{\betabf\ubf})^{-1}(\Kbf_{\ubf\ybf} + \Qbf\Kbf_{\betabf\ybf}) \ybf = \Kbf \ybf
\end{align*}
and the theorem follows.
\end{proof}

As a corollary, we can verify the equivalence of controllers as follows.
\begin{corollary}
The original state-feedback SLS \eqn{SLSdefault} and the realization \eqn{state-feedback-controller_time} in \thm{state-feedback-controller} are equivalent to controller $\Kbf = \Phibf_{\ubf}\Phibf_{\xbf}^{-1}$.
\end{corollary}

\begin{proof}
In the form of \eqn{controller-realization-with-internal-state}, the original SLS realization is
\begin{align*}
\mat{\ubf\\ \betabf} =
\mat{
0 & 0 & z\Phibf_{\ubf}\\
I & 0 & I - z\Phibf_{\xbf}
}
\mat{\xbf\\ \ubf\\ \betabf}
\end{align*}
and the new realization in \eqn{state-feedback-controller_time} is
\begin{align*}
\mat{\ubf\\ \betabf} =
\mat{
0 & 0 & z\Phibf_{\ubf}\\
z^{-1}(zI-A) & -z^{-1}B & 0
}
\mat{\xbf\\ \ubf\\ \betabf}.
\end{align*}
And the equivalence is established by \thm{K-equivalent}.
\end{proof}
The output-feedback versions \eqn{SLS-output-feedback-ybar} and \eqn{output-feedback-controller_time} can be verified in the same way.

\subsection{Deployment Requirements}\label{sec:implementation-deployment-requirements}
When deploying a realization in practice, we have finite storage and computation resources. Any realization that assumes infinite dimension parameters will need to be truncated (by, for example, finite FIR approximation as described in \sec{realization}). Also, we have to be frugal in communication. Unlike a system on a single substrate, information exchange is time-consuming (when sent over geo-distributed networks) and costly (if new links need to be built) in a CPS.
Another subtle constraint brought by a digital cyber layer is that the signals should avoid algebraic loops. Each signal in the system should be computed in a serial order.

We aim to formulate those requirements as constraints in mathematical forms. With \thm{K-equivalent}, we can focus on imposing the constraints on the general realization \eqn{controller-realization-with-internal-state}. Consider a frequency-domain based discrete-time implementation, the controller in \eqn{controller-realization-with-internal-state} can be expressed as
\begin{gather*}
\mat{\Kbf_{\ubf\ybf} & 0 & \Kbf_{\ubf\betabf}\\
\Kbf_{\betabf\ybf} & \Kbf_{\betabf\ubf} & \Kbf_{\betabf\betabf}
}\\
=
\mat{\sum\limits_{\tau\geq 0}^{T_{uy}} z^{-\tau} K_{uy}[\tau] &
0 &
\sum\limits_{\tau\geq 0}^{T_{u\beta}} z^{-\tau} K_{u\beta}[\tau]\\
\sum\limits_{\tau\geq 0}^{T_{\beta y}} z^{-\tau} K_{\beta y}[\tau] &
\sum\limits_{\tau\geq 0}^{T_{\beta u}} z^{-\tau} K_{\beta u}[\tau] &
\sum\limits_{\tau\geq 0}^{T_{\beta \beta}} z^{-\tau} K_{\beta \beta}[\tau]
}
\end{gather*}
where $T_{(\cdot)}$ is the largest $\tau$ such that $K_{(\cdot)}[\tau] \neq 0$.
We can then interpret the requirements as follows.

\paragraph{Finite storage:}
$y[t], u[t], \beta[t]$ are all of finite dimensional signals. When $y[t]$ and $u[t]$ are fixed, the finite storage constraint requires a finite dimensional internal state $\beta[t]$, or $\beta[t] \in \Rmbb^{N_{\beta}}$ for some finite $N_{\beta}$.

\paragraph{Finite computation:}
There are finite non-zero elements in all spectral components in the transfer matrices. Consequently, $T_{uy}, T_{u\beta},T_{\beta y},T_{\beta u},T_{\beta \beta}$ should all be finite.

\paragraph{Limited communications:}
Each dimension in $y[t], u[t], \beta[t]$ can only affect some ``neighboring'' dimensions. In other words, $K_{uy}, K_{u\beta}, K_{\beta y}, K_{\beta u}, K_{\beta \beta}$ are subject to communication constraints, i.e., some elements must be zeros according to the physical structure of the system.

\paragraph{Freedom from algebraic loops for digital cyber layer:}

When computing the signals by modern digital cyber components such as computers, the computations are done sequentially by programs or scripts. If two signals depend on each other, their sequential computations require multiple information exchanges for the values to converge. However, in a CPS, exchanging information through communication is expensive.
As a result, we want to compute each signal in one shot, which requires freedom from algebraic loop among the signals at each time step. It means that at time step $\tau$, if signal $a[\tau]$ is computed from $b[\tau]$, $b[\tau]$ should not depend directly or indirectly on $a[\tau]$.

\vspace*{0.5\baselineskip}

The deployment architectures we propose in the following sections all satisfy the above deployment requirements.

\subsection{Basic Components}\label{sec:implementation-components}
Now we analyze \eqn{controller-realization-with-internal-state} to identify the necessary basic components for CPS controller deployment, which we will assume that the target CPS is able to provide at its nodes.
Firstly, we need storage components to store variables
$y[t], u[t],$ and $\beta[t]$. Also, we need computation components to perform the arithmetic operations described by $K_{(\cdot)}[\tau]$. Lastly, we need communication components to send the variables to where the computation is performed. We summarize those basic components below.

\begin{figure}
\centering
\subcaptionbox{Computation and storage components\label{subfig:architecture-components-computation-and-storage}}{
\includegraphics{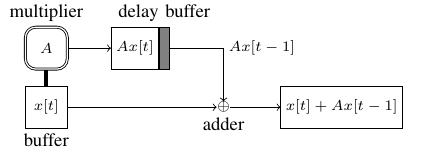}
}\\[0.5\baselineskip]
\subcaptionbox{Communication components\label{subfig:architecture-components-communication}}{
\includegraphics{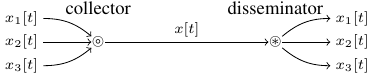}
}
\caption{
Basic components required for deployment. They can be categorized into (\subref{subfig:architecture-components-computation-and-storage}) computation and storage components and (\subref{subfig:architecture-components-communication}) communication components.}
\label{fig:architecture-components}
\end{figure}

A \textbf{\emph{buffer}} is where the system keeps the value of a variable. It can be a memory device such as a register or merely a collection of some buses (wires). In contrast, a \textbf{\emph{delay buffer}} is the physical implementation of $z^{-n}$ in a block diagram. It keeps the variable received at time $t$ and releases it at time $t+n$.

For computation, a \textbf{\emph{multiplier}} takes a vector from a buffer as an input and computes a matrix-vector multiply operation. An \textbf{\emph{adder}} performs entry-wise addition of two vectors of compatible dimensions. Of course, we can merge cascaded adders into a multiple-input adder.

A node can also communicate with other nodes through \emph{disseminator-collector pairs}. A \textbf{\emph{disseminator}} sends some parts of a variable (i.e., a subset of a vector) to designated nodes. At the receiving side, a \textbf{\emph{collector}} assembles the received parts appropriately to reconstruct the desired variable.

To better illustrate the requirements discussed in \sec{implementation-deployment-requirements}, in the following section, we build our controller deployment architecture using the basic components above .

%% file: architecture.tex
\section{Deployment Architectures}\label{sec:architecture}

We now explore the controller architectures for the deployment phase. The crux of our designs centers on the partitions of the controller realization described in \thm{state-feedback-controller} and illustrated in~\fig{simpler-realizations}.
We leverage the basic components described in \sec{implementation-components} and propose the centralized, global state, and distributed architectures accordingly. As long as the real system is capable of providing the basic components, mapping the architectures into the system is straightforward.

%% file: original.tex
\subsection{Architectures for Original SLS Block Diagrams}\label{sec:original}
Internally stabilizing controller realizations are provided in the original SLS papers \cite{AndDLM19,wang2019system}. Using the basic components defined earlier, we demonstrate the architectures corresponding to the original realizations.

\subsubsection{State-Feedback}
\begin{figure}
\centering
\includegraphics{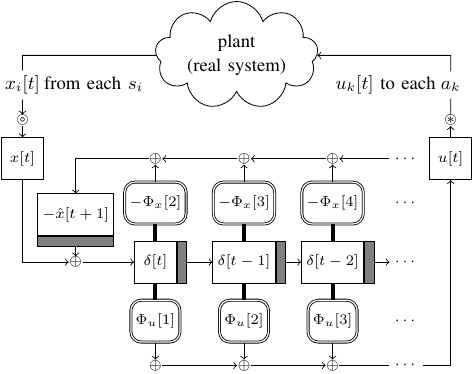}
\caption{The architecture of the original controller model in \fig{SLS-state-feedback-controller-diagram}.
}
\label{fig:original-architecture}
\end{figure}

The original realization of a state-feedback controller is shown in \fig{SLS-state-feedback-controller-diagram}. To derive the corresponding architecture, we consider its state space expression \eqn{SLSdefault} and translate the computation into basic components. The resulting architecture is in \fig{original-architecture}. We note that this controller is applicable to both stable and unstable plants.

\subsubsection{Output-Feedback}
Similarly, we can derive the architecture of the original output-feedback realization in~\fig{SLS-output-feedback-controller-diagram}. A naive state space implementation of \fig{SLS-output-feedback-controller-diagram} in \eqn{SLS-output-feedback-ybar} would lead to an algebraic loop between $u[t]$ and $\overline{y}[t]$. To admit a digital cyber layer, we derive an algebraic-loop-free architecture as follows.
We first separate $u\t$ into two terms by defining
\begin{align*}
u'\t =&\
\sum_{\tau \geq 1}\Phi_{ux}[\tau]\beta[t + 1 - \tau] +
\sum_{\tau \geq 1}\Phi_{uy}[\tau]\overline{y}[t - \tau].
\end{align*}
As a result, we have $u\t =\ \Phi_{uy}[0]\overline{y}[t] + u'\t.$
Substituting \eqn{SLS-output-feedback-ybar} into the above equation, we obtain
\begin{align}\label{eqn:ofbasicrealization}
(I + \Phi_{uy}[0]D)u\t = \Phi_{uy}[0]y\t + u'\t.
\end{align}
Due to lack of space, we omit the full block diagram of this controller. However, it can be intuitively constructed by following the state-feedback case and additionally implementing the left and right hand sides of~\eqn{ofbasicrealization} as a multiplier and buffer respectively.

%% file: centralized.tex
\subsection{Centralized Architecture}\label{sec:centralized}
\begin{figure}
\centering
\includegraphics{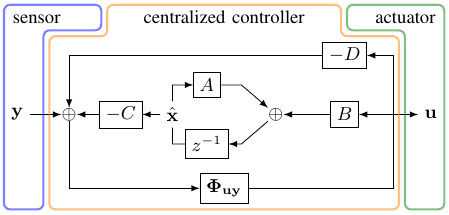}
\caption{The partitions of the output-feedback block diagrams for the centralized architecture. State feedback follows in exactly the same manner -- with fewer components.}
\label{fig:centralized}
\end{figure}

The most straightforward deployment architecture for the new realizations derived in \thm{state-feedback-controller} is the \emph{centralized architecture}, which packs all the control functions into one node, the centralized controller. The partitions of the block diagrams are shown in~\fig{centralized}. The centralized controller maintains communications with all sensors and actuators to collect state/measurement information and dispatch the control signal.

\begin{remark}
In the complexity analysis that follows, we consider worst case problem instances, i.e., we ignore sparsity. In reality, $A$ and $B$ will have some sparsity structure and the SLS program will try to produce sparse/structured spectral components $\Phi_{xx}[\cdot], \Phi_{xy}[\cdot], \Phi_{ux}[\cdot]$, and $\Phi_{uy}[\cdot]$ -- although not all systems can be localized, and localization comes with a performance cost. For systems that are easily localizable, the original SLS controller design (\fig{SLS-state-feedback-controller-diagram} and \fig{original-architecture}) may be more efficient than our analysis shows. But likewise, we do not take sparsity in $(A,B)$ into account which would benefit the new realizations. Our goal is to characterize different implementations for difficult to control systems and arm engineers with the appropriate trade-off information.
\end{remark}

\subsubsection{State-Feedback}
\begin{figure}
\centering
\includegraphics{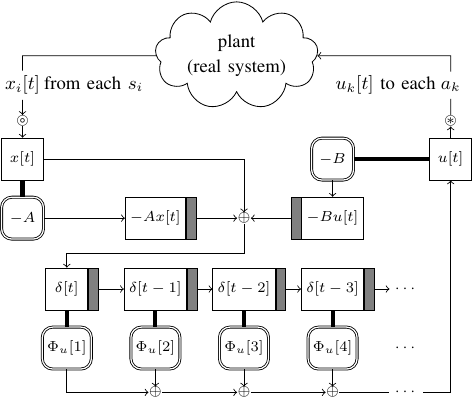}
\caption{The architecture of the centralized state-feedback controller.}
\label{fig:state-feedback-centralized-architecture}
\end{figure}

\fig{state-feedback-centralized-architecture} shows the architecture of the centralized state-feedback controller. For each time step $t$, the controller first collects the state information $x_i\t$ from each sensor $s_i$ for all $i = 1,\dots,N_x$. Along with the stored control signal $u$, the centralized controller computes $\delta\t$ (the disturbance estimate) as in \eqn{state-feedback-state-space-delta}. $\delta\t$ is then fed into an array of delay buffers and multipliers to perform the convolution \eqn{state-feedback-state-space-u} and generate the control signals. The control signals $u_k\t$ are then sent to each actuator $a_k$ for all $k = 1,\dots, N_u$.

Deploying a synthesized solution to the centralized architecture is simple: We take the spectral components ${\Phi_u}[\tau]$ of ${\Phibf_\ubf}$ from $\tau = 1,2,\dots$ and insert them into the array of the multipliers.

We can compare the centralized architecture of \fig{state-feedback-controller-diagram} with the architecture of the original SLS controller (c.f. \fig{original-architecture}).
As mentioned in \sec{background}, the original architecture
is expensive both computationally and storage-wise. Specifically, when both ${\Phibf_\xbf}$ and ${\Phibf_\ubf}$ are FIR with horizon $T$ and not localizable, the original architecture depicted in \fig{original-architecture} performs
\begin{gather*}
\underbrace{\strut (T-1) N_x (2 N_x - 1)}_{\Phi_x[\cdot]\delta[\cdot]} +
\underbrace{\strut T N_u (2N_x - 1)}_{\Phi_u[\cdot]\delta[\cdot]} +
\underbrace{\strut T - 1}_{\text{additions}}
\end{gather*}
floating point operations per time step (flops)\footnote{Our notion of a flop is limited to scalar additions, subtraction, division, and multiplication.
} and needs
\begin{equation}
\underbrace{\strut (T-1)N_x^2}_{\Phi_x[\cdot]} +
\underbrace{\strut T N_x N_u}_{\Phi_u[\cdot]} +
\underbrace{\strut (T + 2) N_x + N_u}_{\delta[\cdot],\ x\t,\ -\hat{x}[t+1], \text{~and~} u\t}\nonumber
\end{equation}
memory locations to store all variables and multipliers. In contrast, the architecture corresponding to \thm{state-feedback-controller} and depicted in \fig{state-feedback-centralized-architecture} requires
\begin{gather}
\underbrace{\strut (N_x + N_u)(2 N_x - 1)}_{-Ax\t\text{~and~}-Bu\t} +
\underbrace{\strut T N_u (2 N_x - 1)}_{\Phi_u[\cdot]\delta[\cdot]} +
\underbrace{\strut T + 1}_{\text{additions}}
\label{eqn:centralized-computation-requirement}
\end{gather}
flops, and needs
\begin{align}
\hspace{-1ex}
\underbrace{\strut N_x^2 + N_x N_u}_{A\text{~and~}B} +
\underbrace{\strut 2 N_x }_{Ax\t\text{,}Bu\t} +
\underbrace{\strut TN_x N_u}_{\Phi_u[\cdot]\delta[\cdot]} +
\underbrace{\strut (T + 1) N_x + N_u}_{\delta[\cdot],\ x\t, \text{~and~} u\t}
\label{eqn:centralized-storage-requirement}
\end{align}
memory locations.

The above analysis tells us that
if the $\Phi_x[\cdot]$ and $\Phi_u[\cdot]$ are dense matrices (i.e., for systems with no localization), the controller architecture corresponding to \thm{state-feedback-controller} is more economic in terms of computation and storage requirements when $N_x \geq N_u$, $N_x \geq 2$, and $T > 3$. We note that having fewer inputs than states is the default setting for distributed control problems. For systems that can be strongly localized (which roughly corresponds to the spectral components having a banded structure with small bandwidth), the original SLS architecture will likely be more economical. A detailed analysis for specific localization regimes is beyond the scope of the paper, here we aim to approximately capture the scaling behavior. It should also be pointed out that we have not taken into account sparsity in $A$ or $B$ for either architectures.

\subsubsection{Output-Feedback}
\begin{figure}
\centering
\includegraphics{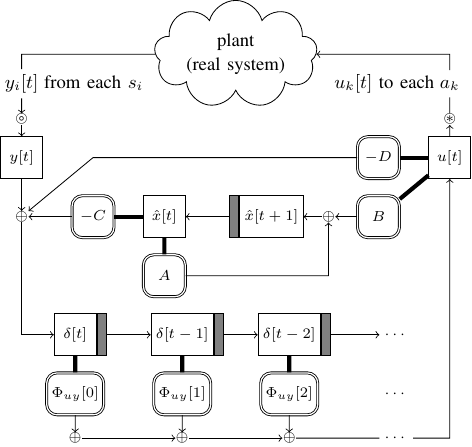}
\caption{The architecture of the centralized output-feedback controller.}
\label{fig:output-feedback-centralized-architecture}
\end{figure}

The centralized output-feedback controller from \thm{state-feedback-controller} is shown in~\fig{output-feedback-centralized-architecture} (due to space limitations, we omit the architecture
diagram corresponding to the standard SLS controller
). The controller collects the measurements $y_i\t$ from each sensor $s_i$ for all $i = 1, \dots, N_y$, maintains internal states $\hat{x}\t$ and $\delta\t$, and generates control signals according to \eqn{output-feedback-controller_time}.

Comparing to the architecture formed by the original SLS controller in~\fig{SLS-output-feedback-controller-diagram},
the new centralized architecture in~\fig{output-feedback-centralized-architecture} is much cheaper in both computation and storage. For FIR $\Phibf_{\xbf\xbf}, \Phibf_{\ubf\xbf}, \Phibf_{\xbf\ybf}$ and $\Phibf_{\ubf\ybf}$ with horizon $T$, the original architecture
requires
\begin{gather*}
(N_u + N_y) (2 N_u - 1) + (T N_u + (T-1) N_x)(2N_x-1) \\+ (N_u + T N_x)(2N_y - 1) + 4T - 1
\end{gather*}
flops and uses
\begin{gather*}
N_u(N_u + 2N_y) + (T-1) (N_u + N_x) (N_x + N_y) \\
+ (T + 1) N_y + T N_x + 2 N_u
\end{gather*}
scalar memory locations, which are more than
\begin{align*}
2(N_x + N_y)(N_x + N_u - 1) + T N_u(2 N_x - 1) + T + 2
\end{align*}
flops and
\begin{gather*}
T N_u N_y + (N_u + N_x) (N_x + N_y) \\
+ (T + 1) N_y + 2 N_x + N_u
\end{gather*}
scalar memory locations used in~\fig{output-feedback-centralized-architecture}. Again, the analysis is based on the most general case. With sparsity or structure, one could potentially further reduce the resource requirement.

\subsubsection{Architecture Trade-Offs}
Despite the intuitive design, the centralized architecture raises several operational concerns. First, the centralized controller becomes the single point of failure in the CPS. Also, the scalability of the centralized scheme is poor: The centralized controller has to ensure communication with all the sensors/actuators and deal with the burden of high computational load. In the following subsections, we explore various distributed system architectures.

%% file: global-state.tex
\subsection{Global State Architecture}\label{sec:global-state}
\begin{figure}
\centering
\includegraphics{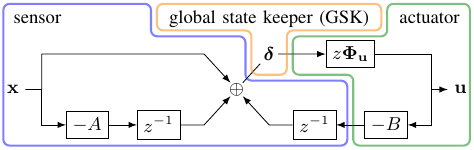}
\caption{The partitions of the state-feedback block diagram for the global state architecture.}
\label{fig:global-state}
\end{figure}

The first attempt to reduce the burden on the controller is to offload the computation to the nodes in the system. A first step in this direction is the \emph{global state architecture}. With this architecture, the sensing and actuation parts of the controller are pushed to the nodes,
and the appropriate information is communicated to a centralized global state keeper (GSK) which keeps track of the full system state instead of the input/output signals at each time $t$. This partitioning is illustrated in \fig{global-state}.

\subsubsection{State-Feedback}
The GSK for a state-feedback controller keeps track of the global state $\delta\t$ instead of the raw state $x\t$. Rather than directly dispatching the control signals $u\t$ to the actuators, GSK supplies $\delta\t$ to the actuators and relies on the actuators to compute $u\t$.

\begin{figure}
\centering
\subcaptionbox{The architecture of the global state keeper (GSK).}{
\includegraphics{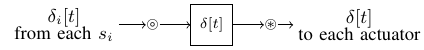}
}\\[0.5\baselineskip]
\subcaptionbox{The architecture of each sensor $i$ ($s_i$).}{
\includegraphics{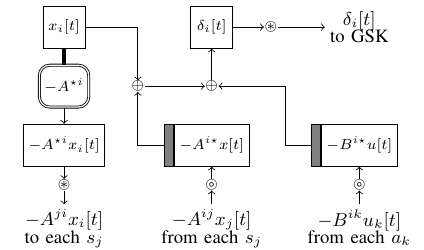}
}\\[0.5\baselineskip]
\subcaptionbox{The architecture of each actuator $k$ ($a_k$).}{
\includegraphics{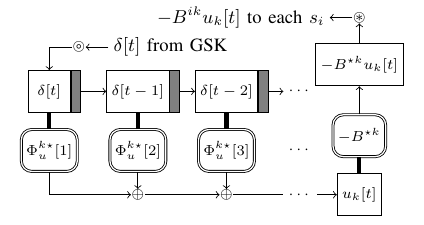}
}
\caption{The state-feedback global state architecture.}
\label{fig:state-feedback-global-state-architecture}
\end{figure}

We illustrate the details of each node in~\fig{state-feedback-global-state-architecture} (and in algorithmic form in~\alg{state-feedback-global-state} where we describe the the computation required to implement a control policy at the sensors, actuators, and GSK). GSK collects $\delta_i\t$ from each sensor $s_i$. To compute $\delta_i\t$, each $s_i$ stores a column vector $-A^{\star i}$. Using the sensed state $x_i\t$, $s_i$ computes $-A^{ji}x_i\t$ and sends it to $s_j$. Meanwhile, $s_i$ collects $-A^{ij}x_j\t$ from each $x_j$ and $-B^{ik}u_k\t$ from each $a_k$. Together, $s_i$ can compute $\delta_i\t$ by
\begin{align*}
\delta_i\t =&\ x_i\t -A^{i\star}x\tm- B^{i\star} u\tm\\
=&\ x_i\t - \sum\limits_{j} A^{ij}x_j\tm - \sum\limits_{k} B^{ik}u_k\tm.
\end{align*}
The $\delta\t$ term is then forwarded to each actuator by GSK. The actuator $a_k$ can compute the control signal using the multiplier array similar to the structure in the centralized architecture. The difference is that each actuator only needs to store the rows of the spectral components ${\Phi_u^{k\star}}[\tau]$ of ${\Phibf_\ubf}$. After getting the control signal $u_k\t$, $a_k$ computes $-B^{ik}u_k\t$ for each sensor $s_i$.

One outcome of this communication pattern is that $s_i$ only needs to receive $-A^{ij}x_j\t$ from $s_j$ if $A^{ij} \neq 0$. Similarly, only when $B^{ik} \neq 0$ does  $s_i$ need to receive information from $a_k$. This property tells us that the nodes only exchange information with their neighbors when a non-zero entry in $A$ and $B$ implies the adjacency of the corresponding nodes (as shown in \subfig{cyber-physical-comparison}{global-state}). Notice that this property holds for any feasible $\SFpair$, regardless of the constraint set $\Scal$.

\begin{algorithm}[!t]
\begin{algorithmic}[1]
{\color{mblue} \algnoindent $\delta[t] = $ \texttt{global\_state\_keeper} $\left(\left\lbrace \delta_i[t] \right\rbrace_i\right)$:}
\STATE{Stack received $\left\lbrace \delta_i[t] \right\rbrace_i$ from all sensors $s_i$
to form $\delta[t]$.}
\end{algorithmic}
\vspace{0.5\baselineskip}

\begin{algorithmic}[1]
{\color{mblue} \algnoindent $\left\lbrace -A^{ji} x_i[t]\right\rbrace_j, \delta_i[t] = $
\item[] \texttt{sensor\_i}$\left(
x_i[t],
\left\lbrace -A^{ij} x_j[t]\right\rbrace_j,
\left\lbrace -B^{ik} u_k[t]\right\rbrace_k
\right)$:}
\STATE{Extract $\left\lbrace -A^{ji} x_i[t]\right\rbrace_j$ from $-A^{\star i} x_i[t]$.}
\STATE{Stack received $\left\lbrace -A^{ij} x_j[t]\right\rbrace_j$ and $\left\lbrace -B^{ik} u_k[t]\right\rbrace_k$ to form $-A^{i\star}x[t]$ and $-B^{i\star}u[t]$.}
\STATE{$\delta_i[t] = x_i[t] - A^{i\star}x[t-1] - B^{i\star}u[t-1]$.}
\end{algorithmic}
\vspace{0.5\baselineskip}

\begin{algorithmic}[1]
{\color{mblue} \algnoindent $\left\lbrace -B^{ik} u_k[t]\right\rbrace_i = $ \texttt{actuator\_k}$(\delta[t])$:}
\STATE{$u_k[t] = \sum\limits_{\tau \geq 1} \Phi_{u}^{k\star}[\tau] \delta[t+1-\tau]$.}
\STATE{Extract $\left\lbrace -B^{ik} u_k[t]\right\rbrace_i$ from $-B^{\star k} u_k[t]$.}
\end{algorithmic}

\caption{The state-feedback global state architecture.}
\label{alg:state-feedback-global-state}
\end{algorithm}

\subsubsection{Architecture Trade-Offs}
The global state architecture mitigates the computation workload of a CPS with a centralized architecture an. However, such an architecture is subject to a single point of failure -- the GSK. We now shift our focus to distributed architectures that are not as vulnerable.

%% file: distributed.tex
\subsection{Distributed Architectures}\label{sec:distributed}
To avoid the single point of failure problem, we can either reinforce the centralized unit through redundancy, or, decompose it into multiple sub-units, each responsible for a smaller region.
The former option trades resource efficiency for robustness, while the latter prevents a total blackout by localizing the failures.
Here we take the latter option to an extreme: We remove GSK from the partitions of the control diagram (\fig{global-state}) and distribute all control functions to the nodes in the network.
As such, the impact of a failed node is contained within a small part of the network.

\subsubsection{State-Feedback}
\begin{figure}
\centering
\includegraphics{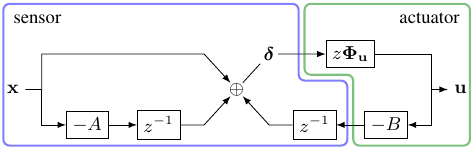}
\caption{A naive way to partition the block diagram in a distributed manner: The sensors directly send $\delta\t$ to each actuator. Those duplicated copies of $\delta\t$ waste memory resources.}
\label{fig:state-feedback-distributed-naive}
\end{figure}

A naive way to remove the GSK from the state-feedback partitions is to send the state $\delta\t$ directly from the sensors to the actuators, as shown in \fig{state-feedback-distributed-naive}. More specifically, each sensor $s_i$ would send $\delta_i\t$ to all the actuators. Each actuator $a_k$ then assembles $\delta\t$ from the received $\delta_j\t$ for all $j$ and computes $u_k\t$ accordingly. We summarize the function of each sensor/actuator in~\alg{state-feedback-distributed-naive}.

\begin{algorithm}[!t]
\begin{algorithmic}[1]
{\color{mblue} \algnoindent $\delta_i[t], \left\lbrace -A^{ji} x_i[t]\right\rbrace_j = $
\item[] \texttt{sensor\_i}$\left(
x_i[t],
\left\lbrace -A^{ij} x_j[t]\right\rbrace_j,
\left\lbrace -B^{ik} u_k[t]\right\rbrace_k
\right)$:}
\STATE{Extract $\left\lbrace -A^{ji} x_i[t]\right\rbrace_j$ from $-A^{\star i} x_i[t]$.}
\STATE{Stack received $\left\lbrace -A^{ij} x_j[t]\right\rbrace_j$ and $\left\lbrace -B^{ik} u_k[t]\right\rbrace_k$ to form $-A^{i\star}x[t]$ and $-B^{i\star}u[t]$.}
\STATE{$\delta_i[t] = x_i[t] - A^{i\star}x[t-1] - B^{i\star}u[t-1]$.}
\end{algorithmic}
\vspace{0.5\baselineskip}

\begin{algorithmic}[1]
{\color{mblue} \algnoindent $\left\lbrace -B^{ik} u_k[t]\right\rbrace_i = $ \texttt{actuator\_k}$\left(\left\lbrace \delta_i[t] \right\rbrace_i\right)$:}
\STATE{Stack received $\left\lbrace \delta_i[t] \right\rbrace_i$ to form $\delta[t]$.}
\STATE{$u_k[t] = \sum\limits_{\tau \geq 1} \Phi_{u}^{k\star}[\tau] \delta[t+1-\tau]$.}
\STATE{Extract $\left\lbrace -B^{ik} u_k[t]\right\rbrace_i$ from $-B^{\star k} u_k[t]$.}
\end{algorithmic}

\caption{The state-feedback naive distributed architecture.}
\label{alg:state-feedback-distributed-naive}
\end{algorithm}

\begin{figure}
\centering
\includegraphics{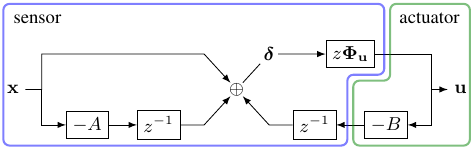}
\caption{The partition of the block diagram for the distributed architecture that conserves memory usage.}
\label{fig:state-feedback-distributed}
\end{figure}

This approach avoids the single point of failure problem. However, it stores duplicated copies of $\delta_i\t$ at each actuator, which wastes memory resources. To conserve  memory usage, we propose to send processed information to the actuators instead of the raw state $\delta\t$. We depict such a memory-conservative distributed scheme in \fig{state-feedback-distributed}.

\begin{figure}
\centering
\subcaptionbox{The architecture of each sensor $i$ ($s_i$).}{
\includegraphics{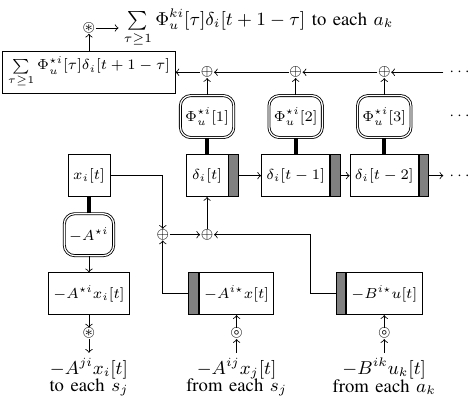}
}\\[0.5\baselineskip]
\subcaptionbox{The architecture of each actuator $k$ ($a_k$).}{
\includegraphics{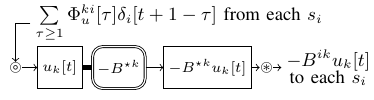}
}
\caption{The memory conservative distributed architecture.}
\label{fig:state-feedback-distributed-architecture}
\end{figure}

The  difference between \fig{state-feedback-distributed-naive} and \fig{state-feedback-distributed} is that we move the convolution $z{\Phibf_\ubf}$ from the actuator side to the sensor side. Implementation-wise, this change leads to the architectures in \fig{state-feedback-distributed-architecture} (and in algorithmic form in~\alg{state-feedback-distributed-memory-conservative}).

In \fig{state-feedback-distributed-architecture}, each sensor $s_i$ not only computes $\delta_i\t$, but $s_i$ also feeds $\delta_i\t$ into a multiplier array for convolution. Each multiplier in the array stores the $i^{\rm th}$ column of a spectral component of ${\Phibf_\ubf}$. The convolution result is then disseminated to each actuator.

At the actuator $a_k$, the control signal $u_k\t$ is given by the sum of the convolution results from each sensor:
\begin{align}
u_k\t =&\ \sum\limits_{\tau \geq 1} {\Phi_u^{k\star}}[\tau]\delta[t+1-\tau] \nonumber \\
=&\ \sum\limits_{i = 1}^{N_x} \sum\limits_{\tau \geq 1} {\Phi_u^{ki}}[\tau]\delta_i[t+1-\tau]. \nonumber
\end{align}

To confirm that \fig{state-feedback-distributed} is more memory efficient than \fig{state-feedback-distributed-naive}, we can count the number of scalar memory locations  in the corresponding architectures. Since both architectures store the matrices $A$, $B$, and $\Phi_u[\cdot]$ in a distributed manner, the total number of scalar memory locations needed for the multipliers is
\begin{align}
N_x^2 + N_x N_u + T N_x N_u.
\label{eqn:distributed-storage-requirement-multiplier}
\end{align}
We now consider the memory requirements for the buffers. Suppose ${\Phibf_\ubf}$ is FIR with horizon $T$ (or, there are $T$ multipliers in the $z{\Phibf_\ubf}$ convolution). In this case \fig{state-feedback-distributed-naive} has the same architecture as \fig{state-feedback-global-state-architecture} without the GSK.
Therefore, the total number of stored scalars (buffers) is as follows:
\begin{align}
\text{At each $s_i$:} &\ N_x + 4, \nonumber\\
\text{At each $a_k$:} &\ (T + 1) N_x + 1, \nonumber\\
\text{Total:} &\ (T + 1) N_x N_u + N_x^2 + 4 N_x + N_u.
\label{eqn:distributed-naive-storage-requirement-buffer}
\end{align}
On the other hand, the memory conservative distributed architecture in \fig{state-feedback-distributed-architecture} uses the following numbers of scalars:
\begin{align}
\text{At each $s_i$:} &\ N_x + N_u + T + 3, \nonumber\\
\text{At each $a_k$:} &\ N_x + 1, \nonumber\\
\text{Total:} &\  2 N_x N_u + N_x^2 + (T + 3) N_x + N_u.
\label{eqn:distributed-storage-requirement-buffer}
\end{align}
In sum, the memory conservative distributed architecture requires fewer scalars in the system, and the difference is
$(T - 1) N_x (N_u - 1)$.

\begin{algorithm}[!t]
\begin{algorithmic}[1]
{\color{mblue} \algnoindent $
\left\lbrace \sum\limits_{\tau \geq 1} \Phi_{u}^{ki}[\tau] \delta_i[t+1-\tau] \right\rbrace_k,
\left\lbrace -A^{ji} x_i[t]\right\rbrace_j = $
\item[] \texttt{sensor\_i}$\left(
x_i[t],
\left\lbrace -A^{ij} x_j[t]\right\rbrace_j,
\left\lbrace -B^{ik} u_k[t]\right\rbrace_k
\right)$:}
\STATE{Extract $\left\lbrace -A^{ji} x_i[t]\right\rbrace_j$ from $-A^{\star i} x_i[t]$.}
\STATE{Stack received $\left\lbrace -A^{ij} x_j[t]\right\rbrace_j$ and $\left\lbrace -B^{ik} u_k[t]\right\rbrace_k$ to form $-A^{i\star}x[t]$ and $-B^{i\star}u[t]$.}
\STATE{$\delta_i[t] = x_i[t] - A^{i\star}x[t-1] - B^{i\star}u[t-1]$.}
\STATE{Extract $\left\lbrace \sum\limits_{\tau \geq 1} \Phi_{u}^{ki}[\tau] \delta_i[t+1-\tau] \right\rbrace_k$ from $\sum\limits_{\tau \geq 1} \Phi_{u}^{\star i}[\tau] \delta_i[t+1-\tau]$.}
\end{algorithmic}
\vspace{0.5\baselineskip}

\begin{algorithmic}[1]
{\color{mblue} \algnoindent $\left\lbrace -B^{ik} u_k[t]\right\rbrace_i = $
\item[] \texttt{actuator\_k}$\left(\left\lbrace \sum\limits_{\tau \geq 1} \Phi_{u}^{ki}[\tau] \delta_i[t+1-\tau] \right\rbrace_i\right)$:}
\STATE{$u_k[t] = \sum\limits_{i=1}^{N_x} \left( \sum\limits_{\tau \geq 1} \Phi_{u}^{ki}[\tau] \delta_i[t+1-\tau] \right)$.}
\STATE{Extract $\left\lbrace -B^{ik} u_k[t]\right\rbrace_i$ from $-B^{\star k} u_k[t]$.}
\end{algorithmic}
\caption{The (state-feedback) memory conservative distributed architecture.}
\label{alg:state-feedback-distributed-memory-conservative}
\end{algorithm}

\subsubsection{Output-Feedback}
\begin{figure}
\centering
\includegraphics{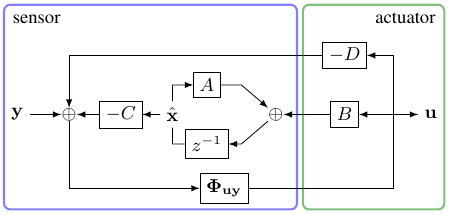}
\includegraphics{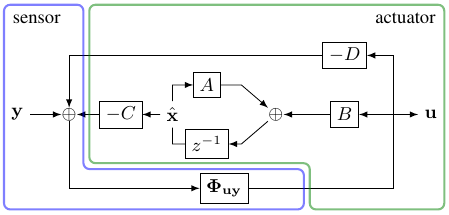}
\caption{Partition choices for the output feedback block diagram for the distributed architecture.}
\label{fig:output-feedback-distributed}
\end{figure}

For output-feedback systems, there are multiple ways to remove the role of a GSK and distribute the computation of the controller. Two possibilities are depicted in~\fig{output-feedback-distributed}. In one, the global state information
is pushed to the sensor, and in the other, to the actuator side.

\begin{figure}
\centering
\includegraphics{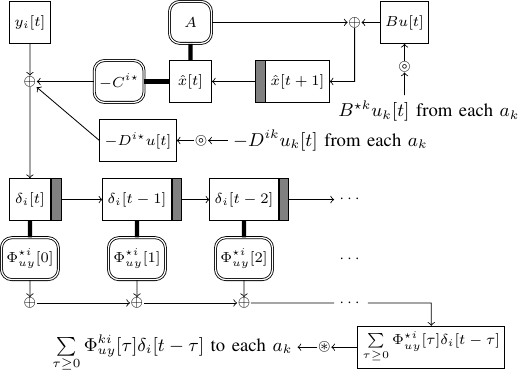}
\caption{The architecture of  sensor $i$ ($s_i$) in the sensor-side global-state distributed architecture. }
\label{fig:output-feedback-distributed-architecture-sensor-side}
\end{figure}

When the sensors take over the global state management, the actuators transmit information
to each sensor, the sensors then build the global state $\hat{x}$ locally as in~\fig{output-feedback-distributed-architecture-sensor-side} (and in algorithmic form in~\alg{output-feedback-distributed-sensor-side}). This approach prevents the single point of failure problem, but it imposes heavy communication load as each sensor  needs to receive two messages; $-D^{ik}u_k\t$, $B^{\star k} u_k\t$ from each actuator $a_k$.

\begin{algorithm}[!t]
\begin{algorithmic}[1]
{\color{mblue} \algnoindent $\left\lbrace \sum\limits_{\tau \geq 0} \Phi_{uy}^{ki}[\tau]\delta_i[t-\tau] \right\rbrace_k  = $
\item[] \texttt{sensor\_i}$\left(
y_i[t],
\left\lbrace B^{\star k} u_k[t]\right\rbrace_k,
\left\lbrace -D^{ik} u_k[t]\right\rbrace_k
\right)$:}
\STATE{Stack received $\left\lbrace B^{\star k} u_k[t]\right\rbrace_k$ to form $Bu[t]$.}
\STATE{Update $\hat{x}[t]$ by \eqn{output-feedback-state-space-hatx}.}
\STATE{Stack received $\left\lbrace -D^{ik} u_k[t]\right\rbrace_k$ to form $-D^{i\star}u[t]$.}
\STATE{$\delta_i[t] = y_i[t] - C^{i\star} \hat{x}[t] - D^{i\star}u[t]$.}
\STATE{Extract $\left\lbrace \sum\limits_{\tau \geq 0} \Phi_{uy}^{ki}[\tau]\delta_i[t-\tau] \right\rbrace_k$ from $\sum\limits_{\tau \geq 0} \Phi_{uy}^{\star i}[\tau]\delta_i[t-\tau]$.}
\end{algorithmic}
\vspace{0.5\baselineskip}

\begin{algorithmic}[1]
{\color{mblue} \algnoindent $\left\lbrace -D^{ik} u_k[t]\right\rbrace_i, B^{\star k} u_k[t] = $
\item[] \texttt{actuator\_k}$\left(\left\lbrace \sum\limits_{\tau \geq 0} \Phi_{uy}^{ki}[\tau]\delta_i[t-\tau] \right\rbrace_i\right)$:}
\STATE{$u_k[t] = \sum\limits_{i = 1}^{N_y} \sum\limits_{\tau \geq 0} \Phi_{uy}^{ki}[\tau]\delta_i[t-\tau]$.}
\STATE{Extract $\left\lbrace -D^{ik} u_k[t]\right\rbrace_i$ from $-D^{\star k} u_k[t]$.}
\STATE{Compute $B^{\star k} u_k[t]$.}
\end{algorithmic}

\caption{The output-feedback sensor-side global-state distributed architecture.}
\label{alg:output-feedback-distributed-sensor-side}
\end{algorithm}

\begin{figure}
\centering
\subcaptionbox{The architecture of  sensor $i$ ($s_i$).}{
\includegraphics{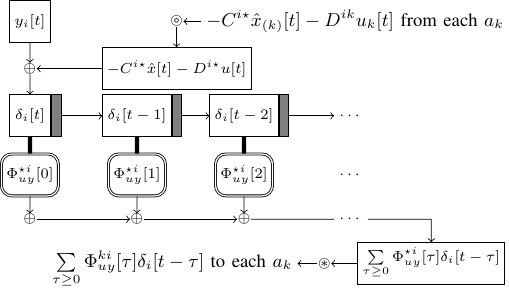}
}\\[0.5\baselineskip]
\subcaptionbox{The architecture of  actuator $k$ ($a_k$).}{
\includegraphics{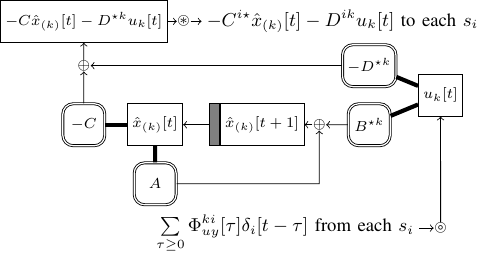}
}
\caption{The actuator-side global-state distributed architecture.}
\label{fig:output-feedback-distributed-architecture-actuator-side}
\end{figure}

One way to alleviate the communication load is to have the actuators send one summarized message instead of two messages to the sensors. However, one copy of $B^{\star k} u_k\t$ is necessary if each sensor needs to compute $\hat{x}\t$, which restrains the actuator from summarizing the information. Therefore, we can merge GSK to the actuator side, which leads to the distributed architecture in~\fig{output-feedback-distributed-architecture-actuator-side} (and in algorithmic form in~\alg{output-feedback-distributed-actuator-side}).

To compute the summarized information, we notice that if we can express $\hat{x}\t$ as the sum of some vectors $\hat{x}_{(k)}\t$ for $k = 1, \dots, N_u$, we have
\begin{align*}
\hat{x}\tn =&\ A\hat{x}\t + Bu\t
= A\sum\limits_{k=1}^{N_u} \hat{x}_{(k)}\t + \sum\limits_{k=1}^{N_u} B^{\star k}u_k\t\\
=&\ \sum\limits_{k=1}^{N_u} \( A\hat{x}_{(k)}\t + B^{\star k}u_k\t \)
\end{align*}
As such, we can enforce $\hat{x}\t = \sum\limits_{k=1}^{N_u} \hat{x}_{(k)}\t$ for all future $t$ by
\begin{align}
\hat{x}_{(k)}\tn = A\hat{x}_{(k)}\t + B^{\star k}u_k\t,
\label{eqn:distributed-actuator-side-hatx}
\end{align}
which defines the internal state $\hat{x}_{(k)}\t$ of each $a_k$. With $\hat{x}_{(k)}\t$, we derive the summarized message at each $a_k$ by
\begin{align}
\delta_i\t =&\ y_i\t - C^{i\star}\hat{x}\t - \sum\limits_{k=1}^{N_u} D^{ik} u_k\t \nonumber \\
=&\ y_i\t + \sum\limits_{k=1}^{N_u} \(- C^{i\star}\hat{x}_{(k)}\t - D^{ik} u_k\t \).
\label{eqn:distributed-actuator-side-delta}
\end{align}

\begin{algorithm}[!t]
\begin{algorithmic}[1]
{\color{mblue} \algnoindent $\left\lbrace \sum\limits_{\tau \geq 0} \Phi_{uy}^{ki}[\tau]\delta_i[t-\tau] \right\rbrace_k  = $
\item[] \texttt{sensor\_$i$}$\left(
y_i[t],
\left\lbrace -C^{i\star}\hat{x}_{(k)}[t] -D^{ik} u_k[t] \right\rbrace_k
\right)$:}
\STATE{Compute $\delta_i[t]$ by \eqn{distributed-actuator-side-delta}.}
\STATE{Extract $\left\lbrace \sum\limits_{\tau \geq 0} \Phi_{uy}^{ki}[\tau]\delta_i[t-\tau] \right\rbrace_k$ from $\sum\limits_{\tau \geq 0} \Phi_{uy}^{\star i}[\tau]\delta_i[t-\tau]$.}
\end{algorithmic}
\vspace{0.5\baselineskip}

\begin{algorithmic}[1]
{\color{mblue} \algnoindent $\left\lbrace -C^{i\star}\hat{x}_{(k)}[t] -D^{ik} u_k[t] \right\rbrace_i = $
\item[] \texttt{actuator\_k}$\left(\left\lbrace \sum\limits_{\tau \geq 0} \Phi_{uy}^{ki}[\tau]\delta_i[t-\tau] \right\rbrace_i\right)$:}
\STATE{$u_k[t] = \sum\limits_{i = 1}^{N_y} \sum\limits_{\tau \geq 0} \Phi_{uy}^{ki}[\tau]\delta_i[t-\tau]$.}
\STATE{Update $\hat{x}_{(k)}[t]$ by \eqn{distributed-actuator-side-hatx}.}
\STATE{Extract $\left\lbrace -C^{i\star}\hat{x}_{(k)}[t] -D^{ik} u_k[t] \right\rbrace_i$ from $-C^{i\star}\hat{x}_{(k)}[t] -D^{ik} u_k[t]$.}
\end{algorithmic}

\caption{The output-feedback actuator-side global-state distributed architecture.}
\label{alg:output-feedback-distributed-actuator-side}
\end{algorithm}

\subsubsection{More Distributed Options}
There exist more distributed architectures under different partitions of the block diagram. For example, one can also assign the output-feedback convolution $\Phibf_{\ubf\ybf}$ to the actuators instead of the sensors. And different partition choices accompany different computation/memory costs. Future work will look at specializing the above computation/memory costs to specific localization constraints. The above expressions serve as upper bounds that are tight for systems that are difficult to localize in space (in the sense of $(d,T)$-localization defined in~\cite{WanMD18}). For systems localizable to smaller regions, the actuators only collect local information, and hence we don't need every sensor to report to every actuator.

%% file: comparison.tex
\section{Architecture Comparison}\label{sec:comparison}

\begin{table*}
\centering
\caption{Comparison Amongst the Proposed Architectures}
\label{tab:comparison}
\def\acol{0.18\textwidth}
\def\bcol{0.175\textwidth}
\renewcommand{\arraystretch}{1.25}

\subcaptionbox{State-Feedback Architectures}{
\begin{tabular}{
|>{\centering}m{\acol}
|>{\centering}m{\bcol}
|>{\centering}m{\bcol}
|>{\centering}m{\bcol}
|>{\centering}m{\bcol}
|}
\hline
\diagbox{Property}{Architecture}
& \textbf{Centralized}
& \textbf{Global State} (\fig{global-state})
& \textbf{Naive Distributed }(\fig{state-feedback-distributed-naive})
& \textbf{Memory Conservative Distributed }(\fig{state-feedback-distributed})
\tabularnewline
\hline
\hline
\textbf{Single Point of Failure}
& yes
& yes
& no
& no
\tabularnewline
\hline
\textbf{Overall Memory Usage}
& lowest\\
(see \eqn{centralized-storage-requirement})
& highest\\
(equal to \eqn{distributed-storage-requirement-multiplier}+\eqn{distributed-naive-storage-requirement-buffer}+$N_x$)
& second highest\\
(equal to \eqn{distributed-storage-requirement-multiplier}+\eqn{distributed-naive-storage-requirement-buffer})
& second lowest\\
(equal to \eqn{distributed-storage-requirement-multiplier}+\eqn{distributed-storage-requirement-buffer})
\tabularnewline
\hline
\textbf{Single Node\\Memory Usage}
& high
& low (actuator needs more memory)
& low (actuator needs more memory)
& low (sensor needs more memory)
\tabularnewline
\hline
\textbf{Single Node\\Computation Loading}
& high (see \eqn{centralized-computation-requirement})
& low (actuator performs convolution)
& low (actuator performs convolution)
& low (sensor performs convolution)
\tabularnewline
\hline
\textbf{Single Node\\Communication Loading}
& sensor/actuator: low\\
controller: high
& sensor/actuator: medium\\
GSK: high
& sensor/actuator:\\high, but localizable
& sensor/actuator:\\high, but localizable
\tabularnewline
\hline
\end{tabular}
}

\vspace*{\baselineskip}

\subcaptionbox{Output-Feedback Architectures (``bot.'' refers to ``bottom'')}{
\begin{tabular}{
|>{\centering}m{\acol}
|>{\centering}m{\bcol}
|>{\centering}m{\bcol}
|>{\centering}m{\bcol}
|>{\centering}m{\bcol}
|}
\hline
\diagbox{Property}{Architecture}
& \textbf{Centralized} (\fig{centralized})
& \textbf{Global State}
& \textbf{Sensor-Side GS Distributed }(\fig{output-feedback-distributed} top)
& \textbf{Actuator-Side GS Distributed }(\fig{output-feedback-distributed} bot.)
\tabularnewline
\hline
\hline
\textbf{Single Point of Failure}
& yes
& yes
& no
& no
\tabularnewline
\hline
\textbf{Overall Memory Usage}
& lowest
& second highest
& highest (sensor stores duplicated $\hat{x}\t$)
& second lowest
\tabularnewline
\hline
\textbf{Single Node\\Memory Usage}
& high
& low
& low (sensor needs more memory)
& low (actuator needs more memory)
\tabularnewline
\hline
\textbf{Single Node\\Computation Loading}
& high
& low
& low (sensor maintains $\hat{x}\t$)
& low (actuator maintains $\hat{x}_{(k)}\t$)
\tabularnewline
\hline
\textbf{Single Node\\Communication Loading}
& sensor/actuator: low\\
controller: high
& sensor/actuator: medium\\
GSK: high
& sensor/actuator:\\very high, but localizable
& sensor/actuator:\\high, but localizable
\tabularnewline
\hline
\end{tabular}
}
\end{table*}

\begin{figure*}
\centering

\subcaptionbox{Centralized architecture.}{
\includegraphics{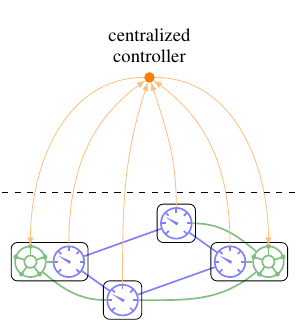}
}\quad\quad
\subcaptionbox{Global state architecture.\label{subfig:cyber-physical-comparison-global-state}}{
\includegraphics{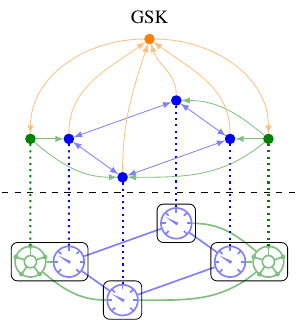}
}\quad\quad
\subcaptionbox{Distributed architecture.}{
\includegraphics{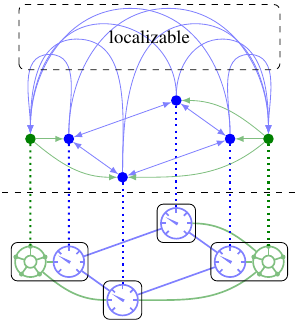}
}

\caption{The cyber-physical structures of the proposed architectures. The horizontal dashed line separates the cyber  and the physical  layers of the system. Each node (the boxes) is equipped with a sensor, which collects the state information in the state-feedback scenarios and the measurements in the output-feedback counterparts.
Solid bullets represents a computation unit, and the arrow links are the communication channels. The cyber structure of the centralized architecture is ignorant about the underlying physical system, while the other architectures manifest some correlations between the cyber and the physical structures. The distributed architectures replace GSK by direct communications, which can be trimmed by imposing appropriate localization constraints on ${\Phibf_\ubf}$ or ${\Phibf_{\ubf\ybf}}$ during the synthesis phase.}
\label{fig:cyber-physical-comparison}
\end{figure*}

We have now seen how different architectures implementing the same controller model allows the engineer to consider different trade-offs when it comes to implementation and deployment. Here we compare the proposed architectures and discuss their differences. Our findings are summarized in \tab{comparison}.

In terms of robustness, the centralized and the global state architectures suffer from a single point of failure, i.e., the loss of the centralized controller or the GSK paralyzes the whole system. This also makes the system vulnerable from a cyber-security perspective. On the contrary, the distributed architectures can still function with some nodes knocked out of the network.\footnote{The system operator might need to update ${\Phibf_\ubf}$ or ${\Phibf_{\ubf\ybf}}$ to maintain  performance --  such a re-design is well within the scope of the SLS framework.}

For information storage, the centralized controller uses the fewest buffers, and the state-feedback global state architecture (because its GSK has to relay $\delta\t$) and the output-feedback sensor-side GS distributed architecture (because the sensors maintain duplicated $\hat{x}\t$) store the most variables. We remark that although the centralized scheme achieves the minimum storage usage at the system level, the single node memory requirement is high for the centralized controller. Conversely, the other architectures store information in a distributed manner, and a small memory is sufficient for each node.

We evaluate the computational load at each node by counting the number of performed multiplication operations. The centralized architecture aggregates all the computation at the centralized controller, while the other architectures perform distributed computing. For distributed settings, the computation overhead is slightly different at each node.
For the state-feedback scenario, the global state and naive distributed architectures let the actuators compute the convolution. Instead, the memory conservative distributed architecture puts the multiplier arrays at the sensors.
For the output-feedback case, convolutions are performed on the sensor side, and the management of global state $\hat{x}\t$ imposes additional computational load on the components. Comparing to the sensor-side GS distributed architecture, it is more computationally balanced for the actuator to handle $\hat{x}\t$.

Finally, we discuss the communication loading of the architectures. In \fig{cyber-physical-comparison}, we sketch the resulting cyber-physical structures of each scheme. In the state-feedback scenarios, each node has its own sensor to measure the state, but some nodes can have no sensors in the output-feedback cases.
The centralized architecture ignores the underlying system interconnection and installs a centralized controller to collect measurements and dispatch control actions. Under this framework, the sensors and the actuators only need to recognize the centralized controller, but the centralized controller must keep track of all the sensors in the system, which limits the scalability of the scheme.

The global state architecture also introduces an additional node into the system, the GSK, with which all  nodes should be contact with. Meanwhile, the sensors and actuators also communicate with each other.
State-feedback sensors and actuators communicate according to the matrices $A$ and $B$. In other words, if two nodes are not directly interacting with each other in the system dynamics, they don't need to establish a direct connection. Output-feedback actuators also communicate with the sensors according to matrix $D$, but sensors compute and distributed control signals to actuators. Although the output-feedback sensors communicate more with the actuators, it is possible to localize the communications by regulating the structure of ${\Phibf_{\ubf\ybf}}$.

Similarly, in the distributed architectures, the sensors and the actuators maintain connections according to $A$ and $B$ for state-feedback or $C$ and $D$ for output-feedback. Additionally, direct communications, which are governed by the structure of ${\Phibf_\ubf}$ or ${\Phibf_{\ubf\ybf}}$, are added to replace the role of GSK. Although it would be slightly more complicated than having a GSK as the relay, we can localize ${\Phibf_\ubf}$ or ${\Phibf_{\ubf\ybf}}$ at the synthesis phase to have a sparse communication pattern.

Besides the centralized architecture, the physical structure has a direct influence on the cyber structure in all other architectures. As such, we believe further research on the deployment architectures of SLS would lead to better cyber-physical control systems.

%% file: future-directions.tex
\section{Conclusion and Future Directions}\label{sec:future}
New internally stabilizing state-feedback and output-feedback controllers were derived for systems that are open-loop stable. The controllers were shown to have block diagram realizations that are in some ways simpler than the ``standard'' SLS controllers.

We considered various architectures to deploy this controller to a real CPS. We illustrated and compared the memory and computation trade-offs among different deployment architectures: centralized, global state, and four different distributed architectures.

There are still many decentralized architecture options left to explore. We are also investigating how  robustness and virtual localizability~\cite{MatWA17} can be integrated into this framework.

%% file: appendices.tex
\section{Proof of \thm{state-feedback-controller}}\label{apx:output-feedback-internal-stability}

The full proof of the state feedback case can be found in~\cite{TseA20}. To provide some insight: observe that  $\Kbf = {\Phibf_\ubf}{\Phibf_\xbf^{-1}} = (z{\Phibf_\ubf})(z{\Phibf_\xbf})^{-1}.$ The realization in \fig{SLS-state-feedback-controller-diagram} then follows by putting $z{\Phibf_\ubf}$ in the forward path and realizing $(z{\Phibf_\xbf})^{-1}$ as the feedback path through the $I-z{\Phibf_\xbf}$ block.

For the output-feedback case, using the Woodbury matrix identity, we know that
\begin{align*}
\Phibf_{\ubf\ybf} =&\ \Lbf + \Lbf C ((zI - A) - B \Lbf C )^{-1} B \Lbf\\
=&\ \(\Lbf^{-1} - C(zI - A)^{-1}B\)^{-1}
= \( \Kbf^{-1} - G \)^{-1}
\end{align*}
where $G = C(zI - A)^{-1}B + D$ is the open-loop transfer matrix. As a result, we know
$\Kbf =\ \( \Phibf_{\ubf\ybf}^{-1} + G \)^{-1}.$
The block diagram then follows immediately.

\begin{figure}
\centering
\includegraphics{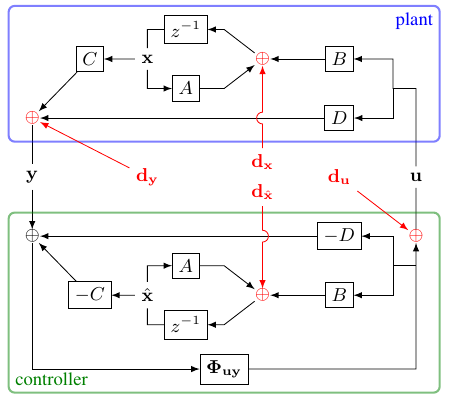}
\caption{Perturbations required for internal stability analysis.
}
\label{fig:output-feedback-internal-stability}
\end{figure}

To show internal stability, we consider additional closed-loop perturbations as in~\fig{output-feedback-internal-stability}. It is sufficient to examine how the internal states $\xbf$, $\ubf$, $\ybf$, and $\hat{\xbf}$ are affected by external signals. The relations among the signals are described by
\begin{align*}
z\xbf =&\ A \xbf + B\ubf + \dbf_\xbf,\\
\ubf =&\ \Phibf_{\ubf\ybf} \(\ybf - C\hat{\xbf} - D (\ubf - \dbf_\ubf)\) + \dbf_\ubf,\\
\ybf =&\ C \xbf + D \ubf + \dbf_\ybf,\\
z\hat{\xbf} =&\ A \hat{\xbf} + B (\ubf - \dbf_\ubf) + \dbf_{\hat{\xbf}}.
\end{align*}
Rearranging the equations above yields the matrix representation
\begin{align*}
\mat{
\xbf \\
\ubf \\
\ybf \\
\hat{\xbf}
} =
\mat{
\Phibf_{\xbf\xbf} & \Delta B \Gamma &
\Phibf_{\xbf\ybf} & \Delta - \Phibf_{\xbf\xbf} \\
\Phibf_{\ubf\xbf} & \Gamma &
\Phibf_{\ubf\ybf} & -\Phibf_{\ubf\xbf}\\
\Lambda C\Delta & \Lambda G &
\Lambda & -G\Phibf_{\ubf\xbf}\\
\Phibf_{\xbf\xbf} - \Delta & \Phibf_{\xbf\ybf}G &
\Phibf_{\xbf\ybf} & 2\Delta - \Phibf_{\xbf\xbf} \\
}
\mat{
\dbf_\xbf \\
\dbf_\ubf \\
\dbf_\ybf \\
\dbf_{\hat{\xbf}}
},
\end{align*}
where
\begin{gather*}
\Delta = (zI-A)^{-1},\quad
\Gamma = \Phibf_{\ubf\ybf}G + I,\quad
\Lambda = G\Phibf_{\ubf\ybf} + I.
\end{gather*}
By assumption, $\Delta$ is stable, and $\Phibf_{\xbf\xbf}, \Phibf_{\ubf\xbf},\Phibf_{\xbf\ybf}, \Phibf_{\ubf\ybf}$ are constrained to be stable by the SLS optimization problem. It can then be verified that all $16$ entries in the matrix are stable, which concludes the proof.

